\documentclass[11pt]{article}

\usepackage{latexsym}
\usepackage{amssymb}

\newtheorem{thm}{Theorem}
\newtheorem{prop}{Proposition}

\newtheorem{dfn}{Definition}

\begin{document}
{
\begin{center}
{\Large\bf
The truncated matrix trigonometric moment problem with an open gap.}
\end{center}
\begin{center}
{\bf S.M. Zagorodnyuk}
\end{center}

\section{Introduction.}
This paper is a continuation of our previous investigations on the truncated matrix trigonometric moment problem (briefly TMTMP)
by the operator approach in~\cite{cit_500_Z},\cite{cit_700_Z}.
The truncated matrix trigonometric moment problem
consists of finding a non-decreasing $\mathbb{C}_{N\times N}$-valued function
$M(t) = (m_{k,l}(t))_{k,l=0}^{N-1}$, $t\in [0,2\pi]$, $M(0)=0$, which is
left-continuous in $(0,2\pi]$, and such that
\begin{equation}
\label{f1_1}
\int_0^{2\pi} e^{int} dM(t) = S_n,\qquad n=0,1,...,d,
\end{equation}
where $\{ S_n \}_{n=0}^d$ is a prescribed sequence of $(N\times N)$ complex matrices (moments).
Here $N\in \mathbb{N}$ and $d\in \mathbb{Z}_+$ are fixed numbers.
Set
\begin{equation}
\label{f1_2}
T_d = (S_{i-j})_{i,j=0}^d =
\left(
\begin{array}{ccccc} S_0 & S_{-1} & S_{-2} & \ldots & S_{-d}\\
S_1 & S_0 & S_{-1} & \ldots & S_{-d+1}\\
S_2 & S_1 & S_0 & \ldots & S_{-d+2}\\
\vdots & \vdots & \vdots & \ddots & \vdots\\
S_d & S_{d-1} & S_{d-2} & \ldots & S_0\end{array}
\right),
\end{equation}
where
$$ S_k := S_{-k}^*,\qquad k=-d,-d+1,...,-1, $$
and $\{ S_n \}_{n=0}^d$ are from~(\ref{f1_1}).
It is well known that the following condition:
\begin{equation}
\label{f1_4}
T_d\geq 0,
\end{equation}
is necessary and sufficient for the solvability of the moment problem~(\ref{f1_1}) (e.g.~\cite{cit_4500_A}).
The  solvable moment problem~(\ref{f1_1}) is said to be \textit{determinate} if it has a unique solution and \textit{indeterminate}
in the opposite case.
We shall omit here an exposition on the history and recent results for the moment problem~(\ref{f1_1}).
All that can be found in~\cite{cit_500_Z},\cite{cit_700_Z}.

Choose an arbitrary $a\in \mathbb{N}$.
Denote by $\mathrm{S}(\mathbb{D};\mathbb{C}_{a\times a})$ a set of all analytic in $\mathbb{D}$, $\mathbb{C}_{a\times a}$-valued functions
$F_\zeta$, such that $F_\zeta^* F_\zeta \leq I_a$, $\forall\zeta\in \mathbb{D}$.
In~\cite{cit_700_Z} we obtained a Nevanlinna-type parameterization for all solutions of the moment problem~(\ref{f1_1}):

\begin{thm}
\label{t2_2}
Let the moment problem~(\ref{f1_1}), with $d\in \mathbb{N}$, be given and condition~(\ref{f1_4}),
with $T_d$ from~(\ref{f1_2}), be satisfied. Suppose that the moment problem is
indeterminate.
All solutions of the moment problem~(\ref{f1_1}) have the following form:
$$ M(t) = (m_{k,l}(t))_{k,l=0}^{N-1},\qquad \int_0^{2\pi} \frac{1}{1-\zeta e^{it}} dM^T (t) $$
\begin{equation}
\label{f1_59}
=
\frac{1}{h_\zeta} \mathbf{A}_{\zeta} - \frac{\zeta}{h_\zeta^2}
\mathbf{B}_{\zeta} F_\zeta
\left(
I_\delta +
\frac{1}{h_\zeta} \mathbf{C}_{\zeta} F_\zeta
\right)^{-1}
\mathbf{D}_{\zeta},\quad \zeta\in \mathbb{D},
\end{equation}
where $\mathbf{A}_{\zeta}$, $\mathbf{B}_{\zeta}$, $\mathbf{C}_{\zeta}$, $\mathbf{D}_{\zeta}$,
are matrix polynomials defined by the given moments, with values in $\mathbb{C}_{N\times N}$, $\mathbb{C}_{N\times \delta}$,
$\mathbb{C}_{\delta\times \delta}$, $\mathbb{C}_{\delta\times N}$, respectively ($0\leq\delta\leq N$).
The scalar polynomial $h_\zeta$ is also defined by the moments.
Here $F_\zeta\in \mathrm{S}(\mathbb{D};\mathbb{C}_{\delta\times \delta})$.
Conversely, each function $F_\zeta\in \mathrm{S}(\mathbb{D};\mathbb{C}_{\delta\times \delta})$
generates by relation~(\ref{f1_59}) a solution of the moment problem~(\ref{f1_1}).
The correspondence between all functions from $\mathrm{S}(\mathbb{D};\mathbb{C}_{\delta\times \delta})$
and all solutions of the moment problem~(\ref{f1_1}) is one-to-one.
\end{thm}

In this paper we shall study the moment problem~(\ref{f1_1}) with an additional constraint posed on the matrix measure $M_{\mathbb{T}}(\delta)$,
$\delta\in \mathfrak{B}(\mathbb{T})$,
generated by the
function $M(x)$ (see the precise definition of $M_{\mathbb{T}}(\delta)$ and other details below):
\begin{equation}
\label{f1_100}
M_{\mathbb{T}}(\Delta) = 0,
\end{equation}
where $\Delta$ is a given open subset of $\mathbb{T}$ (called \textit{a gap}).
Here $\mathbb{T}$ is viewed as a metric space with the metric $r(z,w)=|z-w|$.

\noindent
We present necessary and sufficient conditions for the solvability of the moment problem~(\ref{f1_1}),(\ref{f1_100}). All solutions of
the moment problem~(\ref{f1_1}),(\ref{f1_100}) can be constructed by relation~(\ref{f1_59}), where $F_\zeta$ belongs to a certain subset
of $\mathrm{S}(\mathbb{D};\mathbb{C}_{\delta\times \delta})$.

\textbf{Notations.}
As usual, we denote by $\mathbb{R}, \mathbb{C}, \mathbb{N}, \mathbb{Z}, \mathbb{Z}_+$,
the sets of real numbers, complex numbers, positive integers, integers and non-negative integers,
respectively; $\mathbb{D} = \{ z\in \mathbb{C}:\ |z|<1 \}$; $\mathbb{D}_e = \{ z\in \mathbb{C}:\ |z|>1 \}$;
$\mathbb{T} = \{ z\in \mathbb{C}:\ |z|=1 \}$; $\mathbb{T}_e = \{ z\in \mathbb{C}:\ |z|\not= 1 \}$.
Let $m,n\in \mathbb{N}$.
The set of all complex matrices of size $(m\times n)$ we denote by $\mathbb{C}_{m\times n}$.
The set of all complex non-negative Hermitian matrices of size $(n\times n)$ we denote by $\mathbb{C}_{n\times n}^\geq$.
If $M\in \mathbb{C}_{m\times n}$ then $M^T$ denotes the transpose of $M$, and
$M^*$ denotes the complex conjugate of $M$. The identity matrix from $\mathbb{C}_{n\times n}$
we denote by $I_n$.
By $\mathfrak{B}(\mathbb{T})$ we denote a set of all Borel subsets of $\mathbb{T}$.

If H is a Hilbert space then $(\cdot,\cdot)_H$ and $\| \cdot \|_H$ mean
the scalar product and the norm in $H$, respectively.
Indices may be omitted in obvious cases.
By $\mathbb{C}^N$ we denote the finite-dimensional Hilbert space of complex column vectors
of size $N$ with the usual scalar product $(\vec x,\vec y)_{\mathbb{C}^N} = \sum_{j=0}^{N-1} x_j\overline{y_j}$,
for $\vec x,\vec y\in \mathbb{C}^N$,
$\vec x = (x_0,x_1,\ldots,x_{N-1})^T$, $\vec y = (y_0,y_1,\ldots,y_{N-1})^T$, $x_j,y_j\in \mathbb{C}$.

\noindent For a linear operator $A$ in $H$, we denote by $D(A)$
its  domain, by $R(A)$ its range, by $\mathop{\rm Ker}\nolimits A$
its null subspace (kernel), and $A^*$ means the adjoint operator
if it exists. If $A$ is invertible then $A^{-1}$ means its
inverse. $\overline{A}$ means the closure of the operator, if the
operator is closable. If $A$ is bounded then $\| A \|$ denotes its
norm.
For a set $M\subseteq H$
we denote by $\overline{M}$ the closure of $M$ in the norm of $H$.
For an arbitrary set of elements $\{ x_n \}_{n\in I}$ in
$H$, we denote by $\mathop{\rm Lin}\nolimits\{ x_n \}_{n\in I}$
the set of all linear combinations of elements $x_n$,
and $\mathop{\rm span}\nolimits\{ x_n \}_{n\in I}
:= \overline{ \mathop{\rm Lin}\nolimits\{ x_n \}_{n\in I} }$.
Here $I$ is an arbitrary set of indices.
By $E_H$ we denote the identity operator in $H$, i.e. $E_H x = x$,
$x\in H$. In obvious cases we may omit the index $H$. If $H_1$ is a subspace of $H$, then $P_{H_1} =
P_{H_1}^{H}$ is an operator of the orthogonal projection on $H_1$
in $H$.

By $\mathcal{S}(D;N,N')$ we denote a class of all analytic in a domain $D\subseteq \mathbb{C}$
operator-valued functions $F(z)$, which values are linear non-expanding operators mapping the whole
$N$ into $N'$, where $N$ and $N'$ are some Hilbert spaces.

For a closed isometric operator $V$ in a Hilbert space $H$ we denote:
$M_\zeta(V) = (E_H - \zeta V) D(V)$, $N_\zeta(V) = H\ominus M_\zeta(V)$, $\zeta\in \mathbb{C}$; $M_\infty(V)=R(V)$, $N_\infty(V)= H\ominus R(V)$.

All Hilbert spaces in this paper are assumed to be separable.

\section{The TMTMP with an open gap.}

Let the moment problem~(\ref{f1_1}), with $d\in \mathbb{N}$, be given and condition~(\ref{f1_4}),
with $T_d$ from~(\ref{f1_2}), be satisfied.
Let
$$ T_d = (\gamma_{n,m})_{n,m=0}^{(d+1)N-1},\
S_k = ( S_{k;s,l} )_{s,l=0}^{N-1},\quad -d\leq k\leq d, $$
where $\gamma_{n,m}, S_{k;s,l}\in \mathbb{C}$.
Observe that
\begin{equation}
\label{f2_9}
\gamma_{kN+s,rN+l} = S_{k-r;s,l},\qquad 0\leq k,r\leq d,\quad 0\leq s,l\leq N-1.
\end{equation}
We repeat here some constructions from~\cite{cit_500_Z}.
Consider a complex linear vector space $\mathfrak{H}$, which elements are row vectors
$\vec u = (u_0,u_1,u_2,...,u_{(d+1)N-1})$, with $u_n\in \mathbb{C}$, $0\leq n\leq (d+1)N-1$.
Addition and multiplication by a scalar are defined for vectors in a usual way.
Set
$$ \vec \varepsilon_n = (\delta_{n,0},\delta_{n,1},\delta_{n,2},...,\delta_{n,(d+1)N-1}),\qquad
0\leq n\leq (d+1)N-1, $$
where $\delta_{n,r}$ is Kronecker's delta.
In $\mathfrak{H}$ we define a linear functional $B$ by the following relation:
$$ B(\vec u, \vec w) = \sum_{n,r=0}^{(d+1)N-1} a_n \overline{b_r} \gamma_{n,r}, $$
where
$$ \vec u = \sum_{n=0}^{(d+1)N-1} a_n \vec\varepsilon_n,\quad
\vec w = \sum_{r=0}^{(d+1)N-1} b_r \vec\varepsilon_r,\quad a_n,b_r\in \mathbb{C}. $$
The space $\mathfrak{H}$ with $B$ form a quasi-Hilbert space (\cite{cit_6000_M}).
By the usual procedure of introducing of the classes of equivalence (see, e.g.~\cite{cit_6000_M}), we
put two elements $\vec u$, $\vec w$ from $\mathfrak{H}$ to the same class of equivalence denoted by
$[\vec u]$ or $[\vec w]$, if
$B(\vec u - \vec w,\vec u - \vec w) = 0$. The space of classes of equivalence is a (finite-dimensional) Hilbert
space. Everywhere in what follows it is denoted by $H$.
Set
$$ x_n := [\vec\varepsilon_n ],\qquad 0\leq n\leq (d+1)N-1. $$
Then
\begin{equation}
\label{f2_10}
(x_n,x_m)_H = \gamma_{n,m},\qquad 0\leq n,m\leq (d+1)N-1,
\end{equation}
and $\mathop{\rm span}\nolimits\{ x_n \}_{n=0}^{ (d+1)N-1 } = \mathop{\rm Lin}\nolimits\{ x_n \}_{n=0}^{ (d+1)N-1 }
= H$.
Set $L_N := \mathop{\rm Lin}\nolimits\{ x_n \}_{n=0}^{ N-1 }$.
Consider the following operator:
\begin{equation}
\label{f2_11}
A x = \sum_{k=0}^{dN-1} \alpha_k x_{k+N},\quad x = \sum_{k=0}^{dN-1} \alpha_k x_k,\ \alpha_k\in \mathbb{C}.
\end{equation}

By~\cite[Theorem 1]{cit_500_Z} all solutions of the moment problem~(\ref{f1_1}) have the following form
\begin{equation}
\label{f2_40}
M(t) = (m_{k,j}(t))_{k,j=0}^{N-1},\qquad  m_{k,j}(t) = (\mathbf{E}_t x_k, x_j)_H,
\end{equation}
where $\mathbf{E}_t$ is a left-continuous spectral function of the isometric operator $A$.
Conversely, each left-continuous spectral function of $A$ generates by~(\ref{f2_40}) a solution of the moment problem~(\ref{f1_1}).
The correspondence between all left-continuous spectral functions of $A$ and all solutions of the moment problem~(\ref{f1_1}),
established by relation~(\ref{f2_40}), is one-to-one.

By~\cite[Theorem 3]{cit_500_Z} all solutions of the moment problem~(\ref{f1_1}) have the following form
\begin{equation}
\label{f2_42}
M(t) = (m_{k,j}(t))_{k,j=0}^{N-1},\qquad t\in [0,2\pi],
\end{equation}
where $m_{k,j}$ are obtained from the following relation:
\begin{equation}
\label{f2_43}
\int_0^{2\pi} \frac{1}{1-\zeta e^{it}} dm_{k,j} (t) =
(\left[
E_H - \zeta ( A \oplus \Phi_\zeta )
\right]^{-1} x_k,x_j)_H,\qquad \zeta\in \mathbb{D}.
\end{equation}
Here $\Phi_\zeta\in \mathcal{S}(D;N_0(A),N_\infty(A))$.
Conversely, each $\Phi_\zeta\in \mathcal{S}(D;N_0(A),N_\infty(A))$
generates by relations~(\ref{f2_42})-(\ref{f2_43}) a solution of the moment problem~(\ref{f1_1}).
The correspondence between all $\Phi_\zeta\in \mathcal{S}(D;N_0(A),N_\infty(A))$
and all solutions of the moment problem~(\ref{f1_1}) is one-to-one.

Observe that the right-hand side of~(\ref{f2_43}) may be written as
$(\mathbf{R}_\zeta(A) x_k,x_j)_H$, where $\mathbf{R}_\zeta(A)$ is a generalized resolvent of the isometric operator $A$.
The correspondence between all generalized resolvents of $A$ and all solutions of the moment problem is one-to-one, as well.

Consider an arbitrary solution $M(x)$ of the moment problem~(\ref{f1_1}).
By the construction in~\cite[pp. 791-793]{cit_500_Z}, the corresponding spectral function $\mathbf{E}_t$ in~(\ref{f2_40}) is generated
by the left-continuous orthogonal resolution of unity $\widetilde E_t$ of a unitary operator $\widetilde U_0$ in a Hilbert space
$H_1\supseteq H$. Moreover, the following relation holds:
$$ \widetilde U_0 = U U_0 U^{-1}, $$
where $U$ is a unitary transformation which maps $L^2(M)$ onto $H_1$, and $U_0$ is the operator of multiplication by $e^{it}$ in
$L^2(M)$.

Denote by $\widetilde E (\delta)$, $\delta\in \mathfrak{B}(\mathbb{T})$, the orthogonal spectral measure of the unitary operator
$\widetilde U_0$. The spectral measure $\widetilde E (\delta)$ and the resolution of the identity $\widetilde E_t$ are related
in the following way:
$$ \widetilde E_t = \widetilde E (\delta_t),\qquad \delta_t = \{
z=e^{i\tau}:\ 0\leq\tau < t \},\ t\in [0,2\pi). $$
Therefore the spectral function $\mathbf{E}_t$ and the corresponding spectral measure $\mathbf{E}(\delta)$, $\delta\in \mathfrak{B}(\mathbb{T})$,
satisfy the following relation:
\begin{equation}
\label{f2_100}
\mathbf{E}_t = \mathbf{E} (\delta_t),\qquad  t\in [0,2\pi).
\end{equation}
Relation~(\ref{f2_40}) may be rewritten in the following form:
\begin{equation}
\label{f2_110}
M(t) = ( (\mathbf{E}(\delta_t) x_k, x_j)_H  )_{k,j=0}^{N-1},\qquad t\in [0,2\pi).
\end{equation}
Define the following $\mathbb{C}_{N\times N}^\geq$-valued measure on $\mathfrak{B}(\mathbb{T})$ (i.e. a $\mathbb{C}_{N\times N}^\geq$-valued function on
$\mathfrak{B}(\mathbb{T})$ which is countably additive):
\begin{equation}
\label{f2_120}
M_{ \mathbb{T} }(\delta) = ( (\mathbf{E}(\delta) x_k, x_j)_H  )_{k,j=0}^{N-1},\qquad \delta\in \mathfrak{B}(\mathbb{T}).
\end{equation}
From this definition and relation~(\ref{f2_110}) it follows that
\begin{equation}
\label{f2_125}
M(t) = M_{ \mathbb{T} }(\delta_t),\qquad t\in [0,2\pi).
\end{equation}
Observe that
\begin{equation}
\label{f2_129}
\int_{\mathbb{T}} z^n dM_{ \mathbb{T} } = \int_0^{2\pi} e^{int} dM(t) = S_n,\qquad n=0,1,...,d.
\end{equation}

An arbitrary $\mathbb{C}_{N\times N}^\geq$-valued measure
$\widetilde M_{ \mathbb{T} }(\delta)$, $\delta\in \mathfrak{B}(\mathbb{T})$, satisfying the following relation:
\begin{equation}
\label{f2_135}
\widetilde M_{ \mathbb{T} }(\delta_t) = M(t),\qquad t\in [0,2\pi),
\end{equation}
coincides with the matrix measure $M_{ \mathbb{T} }(\delta)$.
In fact, we may consider the following functions:
$$ f_{k,j}(\delta;\alpha;\widetilde M_{ \mathbb{T} }) = (\widetilde M_{ \mathbb{T} }(\delta) (\vec u_k + \alpha\vec u_j),
\vec u_k + \alpha\vec u_j )_{\mathbb{C}^N} \geq 0, $$
where $\alpha\in \mathbb{C}$, $\delta\in \mathfrak{B}(\mathbb{T})$,
$\vec u_k = (\delta_{k,0},\delta_{k,1},\ldots,\delta_{k,N-1})^T$, $0\leq k,j\leq N-1$.
The scalar measures $f_{k,j}(\delta;\alpha;\widetilde M_{ \mathbb{T} })$ and $f_{k,j}(\delta;\alpha; M_{ \mathbb{T} })$ coincide on $\delta_t$, $t\in [0,2\pi)$.
Therefore they coincide on the minimal generated algebra $Y$, which consists of all finite unions of disjoint sets of the form
$\delta_{t_2, t_1} = \{
z=e^{i\tau}:\ t_2\leq\tau < t_1 \},\ t_1,t_2\in [0,2\pi)$.
Since the Lebesgue continuation is unique, these scalar measures coincide.
On the other hand, the entries of $M_{ \mathbb{T} }(\delta)$ and $\widetilde M_{ \mathbb{T} }(\delta)$ are expressed via $f_{k,j}$ by the
polarization formula. Then $\widetilde M_{ \mathbb{T} }(\delta) = M_{ \mathbb{T} }(\delta)$.

During the investigation of the moment problem~(\ref{f1_1}),(\ref{f1_100}), it is enough to assume that the corresponding moment problem~(\ref{f1_1})
(with the same moments) is \textbf{indeterminate}.
In fact, if the corresponding moment problem~(\ref{f1_1}) has no solutions than the  moment problem~(\ref{f1_1}),(\ref{f1_100}) has no
solutions, as well.
If the corresponding moment problem~(\ref{f1_1}) has a unique solution than this solution can be found explicitly, and then
condition~(\ref{f1_100}) may be verified directly.

\begin{prop}
\label{p3_1}
Let the indeterminate moment problem~(\ref{f1_1}) with $d\in \mathbb{N}$ be given and the operator $A$ in a Hilbert space $H$ be constructed as
in~(\ref{f2_11}). Let $\Delta\in \mathfrak{B}(\mathbb{T})$ be a fixed set. Let $M(x)$, $x\in [0,2\pi]$, be a solution of the moment
problem~(\ref{f1_1}), $M_{ \mathbb{T} }(\delta),\ \delta\in \mathfrak{B}(\mathbb{T})$, be the matrix measure which is defined by~(\ref{f2_120}) with
the corresponding spectral measure $\mathbf{E}(\delta)$, $\delta\in \mathfrak{B}(\mathbb{T})$.
The following two conditions are equivalent:
\begin{itemize}

\item[{\rm (i)}]
$M_{ \mathbb{T} }(\Delta)=0$;

\item[{\rm (ii)}]
$\mathbf{E}(\Delta)=0$.

\end{itemize}
\end{prop}
\textbf{Proof.}
(ii)$\Rightarrow$(i). It follows directly from the definition~(\ref{f2_120}).

\noindent
(i)$\Rightarrow$(ii).
By~(\ref{f2_120}) we may write:
\begin{equation}
\label{f2_140}
M_{ \mathbb{T} }(\Delta) = ( (\mathbf{E}(\Delta) x_k, x_j)_H  )_{k,j=0}^{N-1} = ( ( \widetilde E(\Delta) x_k, x_j)_H  )_{k,j=0}^{N-1} = 0.
\end{equation}
Choose arbitrary numbers $l,m$: $0\leq l,m\leq dN+N-1$. Let
$l = rN+k$, $l = sN+j$, where $0\leq k,j\leq N-1$, $r,s\in \mathbb{Z}_+$.
Then
$$ ( \widetilde E(\Delta) x_l, x_m)_H = ( \widetilde E(\Delta) \widetilde U_0^r  x_k, \widetilde U_0^s x_j)_H
= ( \widetilde U_0^{r-s} \widetilde E(\Delta) x_k, x_j)_H $$
$$ = \int_{ \mathbb{T} } z^{r-s} d( \widetilde E(\delta) \widetilde E(\Delta) x_k, x_j)_H =
\int_{ \mathbb{T} } z^{r-s} d( \widetilde E(\delta\cap\Delta) x_k, x_j)_H = 0. $$
Therefore $\mathbf{E}(\Delta)=0$.
$\Box$

Let $z,w\in \mathbb{T}$: $z\not= w$. Let $z=e^{it}$, $w=e^{iy}$, $0\leq t,y < 2\pi$. If $t<y$, we denote
$$ l(z,w) = \{ u=e^{i\tau}:\ t < \tau < y \}. $$
If $t>y$, we set $l(z,w) = \mathbb{T}\backslash ( l(w,z)\cup \{ w \} \cup \{ z \} )$.
Thus, $l(z,w)$ is an open arc of $\mathbb{T}$ with ends in $z,w$.

Observe that \textit{an arbitrary open subset $\Delta$ of $\mathbb{T}$, $\Delta\not=\mathbb{T}$, is a finite or countable union of disjoint
open arcs of $\mathbb{T}$}.
In fact, suppose that $\zeta_0\in \mathbb{T}$: $\zeta_0\notin\Delta$. Then
$1 \notin \frac{1}{ \zeta_0 } \Delta = \{ u = \frac{1}{ \zeta_0 } w,\ w\in\Delta \}$.
Set
$$ \Omega := \left\{ x = \mathop{\rm Arg}\nolimits z,\ z\in \frac{1}{ \zeta_0 }\Delta \right\} \subseteq (0,2\pi). $$
The set $\Omega$ is an open subset of $\mathbb{R}$. Therefore it is a finite or countable union od disjoint open intervals $l_j\subseteq (0,2\pi)$.
Then $\frac{1}{ \zeta_0 }\Delta$ is a finite or countable union of disjoint
open arcs of $\mathbb{T}$. Consequently, $\Delta$ has the same property.

The following proposition and theorem are simple consequences of Proposition~4.1 and Theorem~4.13 in~\cite{cit_6500_Z}
and the above proved property.

\begin{prop}
\label{p3_2}
Let $V$ be a closed isometric operator in a Hilbert space $H$, and $\mathbf{F}(\delta)$,
$\delta\in \mathfrak{B}(\mathbb{T})$, be its spectral measure.
Let $\Delta$ be an open subset of $\mathbb{T}$, $\Delta\not=\mathbb{T}$.
The following two conditions are equivalent:

\begin{itemize}
\item[(i)]  $\mathbf{F}(\Delta) = 0$;
\item[(ii)] The generalized resolvent $\mathbf{R}_z(V)$, corresponding to the spectral measure
$\mathbf{F}(\delta)$, admits analytic continuation on the set
$\mathbb{D}\cup\mathbb{D}_e\cup \overline{\Delta}$,
where $\overline{\Delta} = \{ z\in \mathbb{C}:\ \overline{z}\in \Delta \}$.
\end{itemize}
\end{prop}

Let $V$  be a closed isometric operator in a Hilbert space $H$, and $\zeta\in \mathbb{T}$. Suppose that
$\zeta^{-1}$ is a point of the regular type of $V$.
Consider the following operators (see~\cite[p. 270]{cit_6500_Z}):
\begin{equation}
\label{f2_142}
W_\zeta P^H_{N_0} f = \zeta^{-1} P^H_{N_\infty} f,\qquad f\in N_\zeta,
\end{equation}
with the domain $D(W_\zeta) = P^H_{N_0} N_\zeta$;
\begin{equation}
\label{f2_144}
S = S_\zeta = P^H_{N_0(V)}|_{N_\zeta(V)},\quad Q = Q_\zeta = P^H_{N_\infty(V)}|_{N_\zeta(V)}.
\end{equation}
Moreover (see~\cite[p. 271]{cit_6500_Z}), since $\zeta^{-1}$ is a point of the regular type of $V$, $S^{-1}$ exists and it is defined on the whole $N_0(V)$,
$D(W_\zeta) = N_0(V)$, and
\begin{equation}
\label{f2_146}
W_\zeta = \zeta^{-1} Q_\zeta S_\zeta^{-1}.
\end{equation}

\begin{thm}
\label{t3_1}
Let $V$    be a closed isometric operator in a Hilbert space $H$, and
$\Delta$ be an open subset of $\mathbb{T}$, $\Delta\not=\mathbb{T}$, such that
\begin{equation}
\label{f2_150}
\zeta^{-1} \mbox{is a point of the regular type of $V$},\qquad \forall\zeta\in\Delta,
\end{equation}
and
\begin{equation}
\label{f2_160}
P^H_{M_\infty(V)} M_\zeta(V) = M_\infty(V),\qquad \forall\zeta\in\Delta.
\end{equation}
Let $\mathbf{R}_z = \mathbf{R}_z(V)$  be an arbitrary generalized resolvent of $V$,
and $F_\zeta\in \mathcal{S}(\mathbb{D}; N_{0},N_\infty)$
corresponds to $\mathbf{R}_z(V)$ by Chumakin's formula.
The operator-valued function $\mathbf{R}_z(V)$ has an analytic continuation
on a set $\mathbb{D}\cup\mathbb{D}_e\cup \Delta$ if and only if
the following conditions hold:

\begin{itemize}
\item[1)] $F_\zeta$  admits a continuation on a set $\mathbb{D}\cup\Delta$ and this continuation
is continuous in the uniform operator topology;

\item[2)] The continued function $F_\zeta$ maps isometrically $N_{0}(V)$ on the whole
$N_\infty(V)$, for all $\zeta\in\Delta$;

\item[3)] The operator $F_\zeta-W_\zeta$ is invertible for all $\zeta\in\Delta$, and
\begin{equation}
\label{f2_170}
(F_\zeta - W_\zeta) N_0(V) = N_\infty(V),\qquad \forall\zeta\in\Delta,
\end{equation}
where $W_\zeta$ is from~(\ref{f2_146}).
\end{itemize}
\end{thm}

As it follows from Remark~4.14 in~\cite[p. 274]{cit_6500_Z}, conditions~(\ref{f2_150}),(\ref{f2_160})
are necessary for the existence of at least one generalized resolvent of $V$,
which admits an analytic continuation on $\mathbb{T}_e\cup\Delta$.

\begin{prop}
\label{p3_3}
Let $V$  be a closed isometric operator in a finite-dimensional Hilbert space $H$, and
$\Delta$ be an open subset of $\mathbb{T}$, $\Delta\not=\mathbb{T}$, such that
condition~(\ref{f2_150}) holds. Then condition~(\ref{f2_160}) holds true.
\end{prop}
\textbf{Proof.} By Corollary~4.7 in~\cite[p. 268]{cit_6500_Z} we may write:
$$ N_\infty(V)\dotplus M_\zeta(V) = H,\qquad \forall\zeta\in\Delta. $$
Applying $P^H_{M_\infty(V)}$ to the both sides of the latter equality we obtain relation~(\ref{f2_160}).
$\Box$

By Proposition~\ref{p3_3} and Theorem~\ref{t3_1} we get the following result.

\begin{thm}
\label{t3_2}
Let $V$  be a closed isometric operator in a finite-dimensional Hilbert space $H$, and
$\Delta$ be an open subset of $\mathbb{T}$, $\Delta\not=\mathbb{T}$, such that
condition~(\ref{f2_150}) holds.
Let $\mathbf{R}_z = \mathbf{R}_z(V)$  be an arbitrary generalized resolvent of $V$,
and $F_\zeta\in \mathcal{S}(\mathbb{D}; N_{0},N_\infty)$
corresponds to $\mathbf{R}_z(V)$ by Chumakin's formula.
The operator-valued function $\mathbf{R}_z(V)$ has an analytic continuation
on a set $\mathbb{D}\cup\mathbb{D}_e\cup \Delta$ if and only if
the following conditions hold:

\begin{itemize}
\item[1)] $F_\zeta$  admits a continuation on a set $\mathbb{D}\cup\Delta$ and this continuation
is continuous in the uniform operator topology;

\item[2)] The continued function $F_\zeta$ maps isometrically $N_{0}(V)$ into
$N_\infty(V)$, for all $\zeta\in\Delta$;

\item[3)] The operator $F_\zeta-W_\zeta$ is invertible for all $\zeta\in\Delta$,
where $W_\zeta$ is from~(\ref{f2_146}).

\end{itemize}
\end{thm}

Let us return to the investigation of the moment problem~(\ref{f1_1}). At first, we shall obtain some necessary conditions for
the solvability of the moment problem~(\ref{f1_1}),(\ref{f1_100}).
We shall use the orthonormal bases constructed in~\cite{cit_700_Z}.

\begin{prop}
\label{p3_4}
Let the indeterminate moment problem~(\ref{f1_1}) with $d\in \mathbb{N}$ be given and the operator $A$ in a Hilbert space $H$ be constructed as
in~(\ref{f2_11}). Let $\Delta$ be an open subset of $\mathbb{T}$, $\Delta\not=\mathbb{T}$.
Let $\mathfrak{A}^\zeta = \{ g_j(\zeta) \}_{j=0}^{\widetilde\tau-1}$, $0\leq \widetilde\tau\leq dN$, be an orthonormal basis in $M_\zeta(A)$,
obtained by the Gram-Schmidt orthogonalization procedure from the following sequence:
$$ x_0 - \zeta x_N, x_{1} - \zeta x_{N+1},..., x_{dN-1} - \zeta x_{dN+N-1}. $$
Here $\zeta\in \overline{\Delta} = \{ z\in \mathbb{C}:\ \overline{z}\in \Delta \}$.
The case $\widetilde\tau=0$ means that $\mathfrak{A}^\zeta = \emptyset$, and $M_\zeta(A) = \{ 0 \}$.
Then the following conditions are equivalent:

\begin{itemize}

\item[{\rm (a)}]
$\zeta \mbox{ is a point of the regular type of $A$},\qquad \forall\zeta\in \Delta$.

\item[{\rm (b)}]
$M_\zeta(A) \not= \{ 0 \}$, and the matrix $\mathcal{M}_{E_H - \zeta A}$ is invertible, for all $\zeta\in \overline{\Delta}$.
Here we denote by $\mathcal{M}_{E_H - \zeta A}$ the matrix of the operator $E_H - \zeta A$
with respect to the bases $\mathfrak{A}_2$, $\mathfrak{A}^\zeta$.
\end{itemize}

Conditions (a), (b) are necessary for the existence of a solution $M(x)$, $x\in [0,2\pi]$,  of the moment
problem~(\ref{f1_1}),  such that $M_{ \mathbb{T} }(\Delta)=0$.
If conditions~(a),(b) are satisfied then $\widetilde\tau = \tau\geq 1$.
\end{prop}
\textbf{Proof.}
The implication $(b)\Rightarrow (a)$ is obvious.
Conversely, suppose that $\Delta$ consists of points of the regular type of $A$.
Choose an arbitrary $\zeta\in \overline{\Delta}$, and set $z_0 = \overline{\zeta}\in \Delta$. Then $(A - z_0 E_H)^{-1} M_\zeta(A) = D(A)$.
If $\mathcal{M}_\zeta(A) = \{ 0 \}$, we would get
$D(A)=\{ 0 \}$, $S_0=0$, and in this case the moment problem would be determinate ($M(x)\equiv 0$). Therefore
$\mathcal{M}_\zeta(A) \not= \{ 0 \}$. The rest is obvious.

Suppose that there exists a solution $M(x)$, $x\in [0,2\pi]$,  of the moment
problem~(\ref{f1_1}),  such that $M_{ \mathbb{T} }(\Delta)=0$, where $M_{ \mathbb{T} }(\delta),\ \delta\in \mathfrak{B}(\mathbb{T})$, is the corresponding
matrix measure.
By Proposition~\ref{p3_1} we get $\mathbf{E}(\Delta)=0$, where $\mathbf{E}(\delta)$, $\delta\in \mathfrak{B}(\mathbb{T})$, is
the corresponding spectral measure.
By Proposition~\ref{p3_2} this means that the corresponding generalized resolvent $\mathbf{R}_z(A)$ admits an analytic continuation on a set
$\mathbb{T}_e\cup \overline{\Delta}$.
In this case, as it was noticed after Theorem~\ref{t3_1} relation~(\ref{f2_150}) holds for $\zeta\in\overline{\Delta}$.
$\Box$

Consider the indeterminate moment problem~(\ref{f1_1}), such as in Proposition~\ref{p3_4}, and suppose that condition~(b) of Proposition~\ref{p3_4}
is satisfied.
Set $\widetilde L := \mathop{\rm Lin}\nolimits\{  g_0(\zeta), g_1(\zeta),..., g_{\tau-1}(\zeta), x_0, x_1,...,x_{N-1} \}$, $\zeta\in\overline{\Delta}$.
Notice that $\widetilde L = H$.
In fact, this follows from the following inclusion, which may be checked by the induction argument:
$\{ x_n \}_{n=0}^{kN+N-1}\subseteq \widetilde L$, $k=0,1,...,d$.

Apply the Gram-Schmidt orthogonalization procedure to the following sequence:
$$ g_0(\zeta), g_1(\zeta),..., g_{\tau-1}(\zeta), x_0, x_1,...,x_{N-1}. $$
Observe that the first $\tau$ elements are already orthonormal.
During the orthogonalization of the rest $N$ elements we shall obtain an orthonormal set $\mathfrak{A}^\zeta_+ = \{ g_j'(\zeta) \}_{j=0}^{\delta-1}$.
Observe that $\mathfrak{A}^\zeta_+$ is an orthonormal basis in $N_\zeta(A)$, $\zeta\in \overline{\Delta}$.

For the operator $A$ in the Hilbert space $H$ and an arbitrary $\zeta\in\overline{\Delta}$, we may construct the
operators $S_\zeta$, $Q_\zeta$ from~(\ref{f2_144}) with $V=A$.
Let $\mathcal{M}_{S_\zeta}$ ($\mathcal{M}_{Q_\zeta}$) be the matrix of the operator $S_\zeta$ ($Q_\zeta$)
with respect to the bases $\mathfrak{A}^\zeta_+$, $\mathfrak{A}_3$ (respectively, to the bases $\mathfrak{A}^\zeta_+$, $\mathfrak{A}_3'$):
$$ \mathcal{M}_{S_\zeta} =
\left(
(S_\zeta g_k'(\zeta), u_j)_H
\right)_{\tau\leq j\leq\tau+\delta - 1,\ 0\leq k\leq\delta-1} $$
\begin{equation}
\label{f2_180}
=
\left(
(g_k'(\zeta), u_j)_H
\right)_{\tau\leq j\leq\tau+\delta - 1,\ 0\leq k\leq\delta-1};
\end{equation}
$$ \mathcal{M}_{Q_\zeta} =
\left(
(Q_\zeta g_k'(\zeta), v_j)_H
\right)_{\tau\leq j\leq\tau+\delta - 1,\ 0\leq k\leq\delta-1} $$
\begin{equation}
\label{f2_190}
=
\left(
(g_k'(\zeta), v_j)_H
\right)_{\tau\leq j\leq\tau+\delta - 1,\ 0\leq k\leq\delta-1},\qquad  \zeta\in \overline{\Delta}.
\end{equation}
Denote by $\widetilde{ W }_\zeta$ the matrix of the operator $W_\zeta$ from~(\ref{f2_146}) with respect to the
bases $\mathfrak{A}_3$, $\mathfrak{A}_3'$. Then
\begin{equation}
\label{f2_200}
\widetilde{ W }_\zeta = \zeta^{-1}
\mathcal{M}_{Q_\zeta} \mathcal{M}_{S_\zeta}^{-1},\qquad  \zeta\in \overline{\Delta}.
\end{equation}

\begin{dfn}
\label{d3_1}
Choose an arbitrary $a\in \mathbb{N}$, $\Delta\subseteq \mathbb{T}$, and let
$Y_\zeta$ be an arbitrary $\mathbb{C}_{a\times a}$-valued function, $\zeta\in \overline{\Delta}$.
By $\mathrm{S} (\mathbb{D}; \mathbb{C}_{a\times a}; \Delta; Y)$ we denote a set of all functions $G_\zeta$ from $\mathrm{S}(\mathbb{D};\mathbb{C}_{a\times a})$
which satisfy the following conditions:
\begin{itemize}

\item[{\rm A)}]
$G_\zeta$ admits a continuation on $\mathbb{D}\cup \overline{\Delta}$, and the continued function $G_\zeta$ is continuous (i.e. each entry of
$G_\zeta$ is continuous);

\item[{\rm B)}]
$G_\zeta^* G_\zeta = I_a$, for all $\zeta\in \overline{\Delta}$;

\item[{\rm C)}]
The matrix $G_\zeta - Y_\zeta$ is invertible for all $\zeta\in \overline{\Delta}$.
\end{itemize}
\end{dfn}

We denote by $\mathcal{S}(\mathbb{D}; N_0(A),N_\infty(A);\Delta; W)$ a set of
all functions $\Phi_\zeta$ from $\mathcal{S}(\mathbb{D}; N_0(A),N_\infty(A))$,
which satisfy conditions~1)-3) of Theorem~\ref{t3_2} with $V=A$ and $\overline{\Delta}$ instead of $\Delta$.

Consider a transformation $\mathbf{T}$ which for an arbitrary function $\Phi_\zeta\in \mathcal{S}(\mathbb{D}; N_0(A),N_\infty(A))$
put into correspondence the following $\mathbb{C}_{\delta\times\delta}$-valued function $F_\zeta$:
\begin{equation}
\label{f2_210}
F_\zeta = \mathbf{T} \Phi =
\left(
(\Phi_\zeta u_k, v_j)_H
\right)_{\tau\leq j,k\leq \tau+\delta-1},\qquad \zeta\in \mathbb{D}.
\end{equation}

The transformation $\mathbf{T}$ is bijective, and it maps $\mathcal{S}(\mathbb{D}; N_0(A),N_\infty(A))$
on the whole $\mathrm{S} (\mathbb{D}; \mathbb{C}_{\delta\times \delta})$.
The following conditions:
$\mathcal{S}(\mathbb{D}; N_0(A),N_\infty(A);\Delta; W) \not=\emptyset$,
and
$\mathrm{S} (\mathbb{D}; \mathbb{C}_{\delta\times \delta}; \Delta; \widetilde W) \not=\emptyset$,
are equivalent.
If $\mathrm{S} (\mathbb{D}; \mathbb{C}_{\delta\times \delta}; \Delta; \widetilde W) \not=\emptyset$, then
\begin{equation}
\label{f2_220}
\mathbf{T} \mathcal{S}(\mathbb{D}; N_0(A),N_\infty(A);\Delta; W) = \mathrm{S} (\mathbb{D}; \mathbb{C}_{\delta\times \delta}; \Delta; \widetilde W).
\end{equation}
All this can be checked directly by the definitions of the corresponding sets.

\begin{thm}
\label{t3_3}
Let the indeterminate moment problem~(\ref{f1_1}) with $d\in \mathbb{N}$ be given and the operator $A$ in a Hilbert space $H$ be constructed as
in~(\ref{f2_11}). Let $\Delta$ be an open subset of $\mathbb{T}$, $\Delta\not=\mathbb{T}$,
and condition~(b)
of Proposition~\ref{p3_4} be satisfied.
The moment problem~(\ref{f1_1}) has a solution $M(x)$, $x\in [0,2\pi]$, such that $M_{ \mathbb{T} }(\Delta)=0$,
if and only if $\mathrm{S} (\mathbb{D}; \mathbb{C}_{\delta\times \delta}; \Delta; \widetilde W)\not=\emptyset$.

If $\mathrm{S} (\mathbb{D}; \mathbb{C}_{\delta\times \delta}; \Delta; \widetilde W)\not=\emptyset$, then formula~(\ref{f1_59})
establishes a one-to-one correspondence between all solutions $M(x)$, $x\in [0,2\pi]$, of the moment problem~(\ref{f1_1}), such that $M_{ \mathbb{T} }(\Delta)=0$,
and all functions $F_\zeta\in \mathrm{S} (\mathbb{D}; \mathbb{C}_{\delta\times \delta}; \Delta; \widetilde W)$.
\end{thm}
\textbf{Proof.}
By Proposition~\ref{p3_4} we obtain that condition~(\ref{f2_150}) holds for $V=A$ and with $\overline{\Delta}$ instead of $\Delta$.

Suppose that the moment problem~(\ref{f1_1}) has a solution $M(x)$, $x\in [0,2\pi]$, such that $M_{ \mathbb{T} }(\Delta)=0$.
As in the proof of Proposition~\ref{p3_4} we conclude that the corresponding generalized resolvent $\mathbf{R}_z(A)$ admits
an analytic continuation on $\mathbb{T}_e\cup \overline{\Delta}$. Let $\Phi_\zeta$ be the function from $\mathcal{S}(\mathbb{D}; N_0(A),N_\infty(A))$
which corresponds to $\mathbf{R}_z(A)$ by Chumakin's formula.
By Theorem~\ref{t3_2} we obtain that $\Phi_\zeta$ belongs to $\mathcal{S}(\mathbb{D}; N_0(A),N_\infty(A);\Delta; W)$.
By~(\ref{f2_220}) we get $\mathrm{S} (\mathbb{D}; \mathbb{C}_{\delta\times \delta}; \Delta; \widetilde W)\not= \emptyset$.

Conversely, suppose that $\mathrm{S} (\mathbb{D}; \mathbb{C}_{\delta\times \delta}; \Delta; \widetilde W)\not= \emptyset$.
By~(\ref{f2_220}) we can choose a function $\Phi_\zeta\in\mathcal{S}(\mathbb{D}; N_0(A),N_\infty(A);\Delta; W)$.
Let $\mathbf{R}_z(A)$ be the generalized resolvent of $A$ which corresponds to $\Phi_\zeta$ by Chumakin's formula.
By Theorem~\ref{t3_2} we obtain that $\mathbf{R}_z(A)$ admits
an analytic continuation on $\mathbb{T}_e\cup \overline{\Delta}$.
By Proposition~\ref{p3_2} we conclude that $\mathbf{E}(\Delta)=0$.
Finally, applying Proposition~\ref{p3_1} we obtain that $M_{ \mathbb{T} }(\Delta)=0$.

Consider the case $\mathrm{S} (\mathbb{D}; \mathbb{C}_{\delta\times \delta}; \Delta; \widetilde W)\not= \emptyset$.

Choose an arbitrary function $F_\zeta\in\mathrm{S} (\mathbb{D}; \mathbb{C}_{\delta\times \delta}; \Delta; \widetilde W)$.
Let $\Phi_\zeta\in\mathcal{S}(\mathbb{D}; N_0(A),N_\infty(A);\Delta; W)$ be such that $\mathbf{T}\Phi_\zeta = F_\zeta$.
Repeating the above arguments we conclude that for the corresponding solution $M(x)$, $x\in[0,2\pi]$, it holds
$M_{ \mathbb{T} }(\Delta)=0$. By the construction of formula~(\ref{f1_59}), for $F_\zeta$ it corresponds namely $M(x)$.

Conversely, choose an arbitrary solution $M(x)$, $x\in[0,2\pi]$, of the moment problem~(\ref{f1_1}) such that
$M_{ \mathbb{T} }(\Delta)=0$.
Repeating the arguments at the beginning of the proof we obtain that
$\Phi_\zeta$ belongs to $\mathcal{S}(\mathbb{D}; N_0(A),N_\infty(A);\Delta; W)$, where $\Phi_\zeta$ is related to the
corresponding generalized resolvent by Chumakin's formula. Then
$F_\zeta = \mathbf{T}\Phi_\zeta\in \mathrm{S} (\mathbb{D}; \mathbb{C}_{\delta\times \delta}; \Delta; \widetilde W)$.
Observe that $F_\zeta$ corresponds to $M(x)$ by formula~(\ref{f1_59}).

The correspondence between solutions and functions from $\mathrm{S} (\mathbb{D}; \mathbb{C}_{\delta\times \delta}; \Delta; \widetilde W)$
is one-to-one, as it follows from Theorem~\ref{t2_2}.
$\Box$

\noindent
{\bf Example 2.1. } Let $N=3$, $d=1$,
$S_0 = \left(
\begin{array}{ccc} 1 & 1 & 0\\
1 & 1 & 0\\
0 & 0 & 1\end{array}
\right)$,
$S_1 = \left(
\begin{array}{ccc} 1 & 1 & 0\\
1 & 1 & 0\\
0 & 0 & 0\end{array}
\right)$.
The moment problem~(\ref{f1_1}) with moments $S_0,S_1$ was studied in~\cite{cit_700_Z}.
We shall use orthonormal bases and other objects constructed therein.
The following formula:
\begin{equation}
\label{f2_230}
\int_0^{2\pi} \frac{1}{1-\zeta e^{it}} dM^T (t)
= \left(\begin{array}{ccc} \frac{1}{1-\zeta} & \frac{1}{1-\zeta} & 0 \\
\frac{1}{1-\zeta} & \frac{1}{1-\zeta} & 0 \\
0 & 0 & 1 + \zeta^2 \frac{F_\zeta}{1-\zeta^2 F_\zeta}\end{array}\right),\ \zeta\in \mathbb{D},
\end{equation}
establishes a one-to-one correspondence between all solutions $M(x)$, $x\in [0,2\pi]$, of the moment problem~(\ref{f1_1})
and all functions $F_\zeta\in \mathrm{S} (\mathbb{D}; \mathbb{C}_{1\times 1})$.

\noindent
Let $\Delta = l(1,-1)$. Let us find solutions of the moment problem~(\ref{f1_1}) which satisfy condition~(\ref{f1_100}).
Calculate the elements of the orthonormal bases $\mathfrak{A}^\zeta$ and $\mathfrak{A}^\zeta_+$:
$$ g_0 = \frac{1}{ \sqrt{2-\overline{\zeta}-\zeta} } (x_0 - \zeta x_3),\quad
g_1 = \frac{1}{ \sqrt{2} } (x_2 - \zeta x_5); $$
$$ g_1 = \frac{1}{ \sqrt{2} } (x_2 + \zeta x_5),\qquad \zeta\in\overline{l(1,-1)}. $$
By~(\ref{f2_180}),(\ref{f2_190}) we get:
$$ \mathcal{M}_{Q_\zeta} = \frac{1}{\sqrt{2}},\ \mathcal{M}_{S_\zeta} = \frac{1}{\sqrt{2}} \zeta,\qquad \zeta\in\overline{l(1,-1)}. $$
By~(\ref{f2_200}) we get:
$$ \widetilde W_\zeta = \zeta^{-2},\qquad \zeta\in\overline{l(1,-1)}. $$
Observe that the function $F_\zeta\equiv 1$ belongs to the set $\mathrm{S} (\mathbb{D}; \mathbb{C}_{1\times 1}; l(1,-1); \zeta^{-2})$.
By Theorem~\ref{t3_3} the moment problem~(\ref{f1_1}) with an additional constraint~(\ref{f1_100}) is solvable.
Formula~(\ref{f2_230}) establishes a one-to-one correspondence between all solutions $M(x)$, $x\in [0,2\pi]$, of the moment problem~(\ref{f1_1})
with an additional constraint~(\ref{f1_100})
and all functions $F_\zeta\in \mathrm{S} (\mathbb{D}; \mathbb{C}_{1\times 1}; l(1,-1); \zeta^{-2})$.

\begin{center}
{\large\bf The truncated matrix trigonometric moment problem with an open gap.}
\end{center}
\begin{center}
{\bf S.M. Zagorodnyuk}
\end{center}

This paper is a continuation of our previous investigations on the truncated matrix trigonometric moment problem
in 
Ukrainian Math. J., 2011, \textbf{63}, no. 6, 786-797, 
and 
Ukrainian Math. J., 2013, \textbf{64}, no. 8, 1199-1214.
In this paper we shall study the truncated matrix trigonometric moment problem with an additional constraint posed on the
matrix measure $M_{\mathbb{T}}(\delta)$,
$\delta\in \mathfrak{B}(\mathbb{T})$,
generated by the seeked
function $M(x)$: $M_{\mathbb{T}}(\Delta) = 0$, where $\Delta$ is a given open subset of $\mathbb{T}$ (called a gap).
We present necessary and sufficient conditions for the solvability of the moment problem with a gap. All solutions of
the moment problem with a gap can be constructed by a Nevanlinna-type formula.

}


\end{document}